% Logic Eprints
%Submitted 1641 Fri Mar 10, 1995 by: jlh@cs.cs.appstate.edu (jeff hirst)
%logic/hirstlempp/infver.tex
%

%Article:  Infinite versions of some problems from
%          finite complexity theory
%TeX Format:  AMSTeX
%Authors:  J. Hirst and S. Lempp
%Submitted to:  Illinois Journal of Mathematics
%Author's email:  jlh@math.appstate.edu
\input amstex
\documentstyle{amsppt}
%
%	modifications to amsppt style
% 
\magnification=\magstep1  	%This modifies the size of the print
\hsize = 6.5 truein      	%This modifies page size
\vsize = 9 truein         	%This modifies page size
\TagsAsMath              	%Modification to tags (see Joy, p186)
\NoRunningHeads			%Omit silly running page headers
%	For other modifications, see the Special Definitions section
%	following the \document command, and the double spacing line
%	just prior to the \document command.
%
\topmatter
\title Infinite Versions of Some Problems \\ From Finite Complexity Theory
\endtitle
\author  
Jeffry L. Hirst\\
{\rm Appalachian State University}\\
\qquad\\
Steffen Lempp\\
{\rm University of Wisconsin}
\endauthor
\address
Dept. of Mathematical Sciences,
Appalachian State University,
Boone, NC  28608
\endaddress
\email
jlh\@math.appstate.edu
\endemail
\address
Dept. of Mathematics
University of Wisconsin
Madison, WI  53706-1388
\endaddress
\email
lempp\@math.wisc.edu
\endemail
\date
February 23, 1995
%version 2/23 noon
\enddate
\thanks
Lempp's research was partially supported by NSF grant DMS-9100114.
\endthanks
\endkeywords
\endsubjclass
\abstract
Recently, several authors have explored the connections between NP-complete
problems for finite objects and the complexity of their analogs for infinite objects.
In this paper, we will categorize infinite versions of several problems arising
from finite complexity theory in terms of their recursion theoretic complexity and
proof theoretic strength.  These infinite analogs can behave in a variety of
unexpected ways.
\endabstract
\endtopmatter
%
%	Double space body of paper
\addto\tenpoint{\normalbaselineskip20pt\normalbaselines} 	% extra space between lines
\document
%Special definitions used in this paper.
\define\rca{\bold {RCA_0}}

\define\aca{\bold {ACA_0}}
\define\s11{\bold{\Sigma_1^1{-}AC_0}}

\define\p11{\bold {\Pi_1^1{-}CA_0}}

Startling parallels exist between the computational complexity of certain
graph theoretic problems and the recursion theoretic complexity and
proof theoretic strength of their infinite analogs.  For example,
the problem of deciding which finite graphs have an Euler path is known
to be P-time computable \cite{9}, and {\smc Beigel} and {\smc Gasarch}
\cite{4} have
shown that the problem of deciding which infinite recursive graphs have
an Euler path is arithmetical.  By contrast, the problem of deciding
which finite graphs have Hamilton paths is NP-complete \cite{8},
and {\smc Harel} \cite{6} has shown that the problem of deciding which
infinite recursive graphs have a Hamilton graph is $\Sigma_1^1$ complete.
Thus, the possibly greater computational complexity is paralleled
by a demonstrable increase in recursion theoretic complexity.
This pattern can also be seen through an application of the techniques
of reverse mathematics.  The existence of a function that decides which
graphs have Euler paths is provably equivalent to $\aca$, while the
existence of a similar function for Hamilton paths is equivalent
to the much stronger axiom system $\p11$.

Unfortunately, other graph theoretic problems do not demonstrate this
parallelism.  We have selected some examples to illustrate two
general themes.  First, {\sl different infinite formulations of a fixed
finite problem can have different recursion theoretic complexities.}
This would seem to indicate that the use of a preferred infinite
formulation might lead to natural parallels between finite complexity
and recursion theoretic complexity.  However, the behavior of infinite
analogs is not so easily tamed.  Indeed, {\sl similar formulations of infinite
versions of problems with different finite complexities may have the same
recursion theoretic complexity.}

\subhead{Variability due to translations}\endsubhead

This section contains examples illustrating our first theme.  The
problem of determining which finite graphs are 3-chromatic is
NP-complete \cite{8}.  Extrapolating from the problem of finding
Hamilton paths, we would expect infinite analogs of the 3-coloring
problem to be $\Sigma_1^1$ complete.  However, the actual recursion
theoretic complexity depends on the formulation of the infinite
analog, as demonstrated by the following three theorems.  Our notation
is patterned after that of {\smc Soare} \cite{12}.

\proclaim{Theorem 1}
{\rm ({\smc Beigel} and {\smc Gasarch} \cite{1})}
The set of indices of 3-chromatic recursive graphs is
$\Pi^0_1$ complete.
\endproclaim

\demo{Proof}
Let $\Cal G_1$ denote the set of indices of 3-chromatic recursive graphs.
Note that $x \in \Cal G_1$ if and only if every finite subgraph of the graph
with index $x$ is 3-chromatic.  Thus, $\Cal G_1$ is a $\Pi_1^0$ definable
subset of the set of indices of recursive graphs.

To show that $\Cal G_1$ is $\Pi_1^0$ complete, let $\Cal G_0$ denote
the set of indices of recursive graphs which are not 3-chromatic.
It suffices to show that $( K, \bar K ) \le_1 (\Cal G_0 , \Cal G_1 )$.
Here $K= \{ e : e \in W_e \}$ is the self-halting set.

For each $e \in \omega$, define the graph $G_e$ as follows.  The vertex set
of $G_e$ is $\omega$, and for $m<n$, the edge $(m,n)$ is in $G_e$ if and
only if $\{e\}(e)$ halts by stage $m$.  For every $e$, $G_e$ is recursive.
By the {\it s-m-n} Theorem, there is a 1-1 recursive function $f$ such that
for every $e$, $f(e)$ is an index for $G_e$.

Note that if $e \in K$, $G_e$ contains an infinite clique.  In this case,
$G_e$ is not 3-chromatic, so $f(e) \in \Cal G_0$.  On the other hand,
if $e \in \bar K$, $G_e$ has no edges.  Such a graph is certainly 3-chromatic,
so $f(e) \in \Cal G_1$.  Thus, $f$ witnesses that
$( K, \bar K ) \le_1 (\Cal G_0 , \Cal G_1 )$, as desired.  \qed
\enddemo

\proclaim{Theorem 2}
The set of indices of recursive graphs with finitely colorable connected
components is $\Pi^0_3$ complete.
\endproclaim

\demo{Proof}
Let $\Cal G_1$ denote the set of indices of recursive
graphs with finitely colorable connected components.
Suppose that $x$ is the index of a graph $G$.  Then
$x \in \Cal G_1$ if and only if for every vertex $v$ of $G$,
there is an integer $k$ such that every finite connected subgraph
of $G$ containing $v$ is $k$-chromatic.  Thus, $\Cal G_1$ is a
$\Pi^0_3$ definable subset of the set of indices of recursive graphs.

To show that $\Cal G_1$ is $\Pi_3^0$ complete,
let $\Cal G_0$ denote the set of indices of those recursive graphs
which have connected components that are not finitely colorable.
It suffices to show that
$( Cof , \overline {Cof} )\le_1 (\Cal G_0 , \Cal G_1 )$.
Here, $Cof = \{ e : W_e \,\text{is cofinite} \}$.

For each $e \in \omega$, define the graph $G_e$ as follows.
$G_e$ will contain vertices labeled $v_{m,n}$ for each $m$ and $n$ in $\omega$,
and some additional unlabeled vertices.
For each $m$, the vertex $v_{m,0}$ will be included in a complete
graph on $m+1$ vertices.  For every $m$ and $j$ all edges of the form
$( v_{m,j} , v_{m,j+1} )$ will be included in $G_e$.  Finally, the edge
$(v_{m,j} , v_{m+1,j} )$ will be included in $G_e$ if and only if
$\{e\}(m)$ halts by stage $j$.  For every $e$, $G_e$ is recursive.
By the {\it s-m-n} Theorem, there is a 1-1 recursive function $f$
such that for every $e$, $f(e)$ is an index for $G_e$.

Note that if $e \in  Cof$, then there is a $j$ such that the vertices
$\{ v_{m,n} : m>j \}$ are all in the same connected component.
Consequently, arbitrarily large complete finite subgraphs are
contained in this component, and it is not finitely colorable.
Thus, if $e \in Cof$, $f(e) \in \Cal G_0$.
Now suppose that $e \in \overline {Cof}$ and $C$ is a connected component of
$G_e$.  $C$ must contain a vertex of the form $v_{m,0}$.
Since $e \in \overline {Cof}$,
there is a least $j$ greater than $m$ such that
$\{e\}(j)$ never halts.  Consequently, $C$ cannot contain any vertex
$v_{n,k}$ such that $n>j$.  This ensures that $C$ is $j+1$-chromatic,
so $f(e) \in \Cal G_1$.  Thus, $f$ witnesses that
$( Cof , \overline {Cof} )\le_1 (\Cal G_0 , \Cal G_1 )$, as desired.
\qed
\enddemo

For the next proof,
we will need the following notation for finite sequences of natural numbers.
Assuming a recursive bijection between $\omega$ and $\omega^{< \omega}$,
we will use a Greek letter (usually $\sigma$ or $\tau$) to denote
both a sequence and its integer code.  The formula $\sigma \subseteq \tau$ means
that $\sigma$ is a (not necessarily proper) initial segment of $\tau$.
Thus, $T$ is a {\it tree} if whenever $\tau \in T$ and $\sigma \subseteq \tau$,
then $\sigma \in T$.

Given an arbitrary index $e$, $\{ e \}$ may or may not be the characteristic
function for a recursive tree.  To streamline our discussion, consider
the following auxiliary function.

\definition{Definition 3}
For $e \in \omega$, the partial recursive function
$\eta_e$ is defined by:
$$
\eta_e (\tau ) =
\cases
1 &\text{if \quad} \forall \sigma \subseteq \tau \, (\{ e \} ( \sigma )= 1),\\
0 &\text{if \quad} 
\{e\}(\tau) = 0 \land
[ \forall \sigma \subseteq \tau \, (\{e\}(\sigma )=0 \lor \{e\}( \sigma ) = 1 ]\\
&\qquad \land [\forall \sigma \subseteq \tau \, \forall \alpha \subseteq \sigma \,
(\{e\}(\alpha ) = 0 \rightarrow \{e\}(\sigma ) = 0 )],\\
\uparrow &\text{otherwise.}
\endcases
$$
\enddefinition
Na\"\i vely, $\eta _e$ approximates the characteristic function of a tree.
In particular, $\eta _e$ is total if and only if $e$ is the index of a recursive
tree.  Note that by the {\it s-m-n} Theorem, there is a 1-1 recursive function
which maps each $e$ to an index for $\eta _e$.

\proclaim{Theorem 4}
{\rm ({\smc Tirza Hirst} and {\smc Harel} \cite{7})}
The set of indices of recursive graphs with colorings which use one
color infinitely often is $\Sigma^1_1$ complete.
\endproclaim

\demo{Proof}
Let $\Cal G$ denote the set of indices of recursive graphs with
colorings that use one color infinitely often.  Note that $x \in \Cal G$
if and only if there is a function $\chi$ mapping the vertices of the graph
with index $x$ into $\omega$, such that $\chi$ maps neighboring vertices to
different values, and $0$ appears infinitely often in the range of $\chi$.
This statement can be formalized using a single existential set quantifier
followed by an arithmetical formula, so $\Cal G$ is $\Sigma_1^1$ definable.

To show that $\Cal G$ is $\Sigma_1^1$ complete, we will show that
$\Cal T \le_1 \Cal G$, where $\Cal T$ denotes the set of indices of
recursive trees which are not well-founded.  With each $e \in \omega$,
we associate a partial recursive graph, $G_e$.  The vertex set for $G_e$ consists
of (codes for) elements of $\omega^{< \omega}$.  For every
$\sigma , \tau \in \omega^{ < \omega }$, the characteristic function for
the edge set of $G_e$ is defined by
$$
E_e ( \sigma , \tau ) =
\cases
0 & \text{ if \quad} \eta_e (\sigma ) = 1 \land \eta_e (\tau ) =1 \land
    (\sigma \subseteq \tau \lor \tau \subseteq \sigma), \\
1 & \text{ if \quad} \eta_e (\sigma ) = 1 \land \eta_e (\tau ) =1 \land
    \neg (\sigma \subseteq \tau \lor \tau \subseteq \sigma), \\
1 & \text{ if \quad} \eta_e (\sigma ) \downarrow \land \eta_e(\tau)\downarrow \land
    (\eta_e(\sigma) = 0 \lor \eta_e(\tau)=0),\\
\uparrow & \text{ otherwise.}
\endcases
$$
Roughly, we connect $\sigma$ and $\tau$ by an edge if they are incomparable
nodes on the tree or if one of them is not in the tree, ignoring those nodes
whose status is suspect.  By the {\it s-m-n} Theorem, there is a 1-1 recursive
function $f$ such that for every $e$, $f(e)$ is an index for $G_e$.

If $e \in \Cal T$, then $e$ is the index of a recursive tree containing an
infinite path $P$.  Consequently, $f(e)$ is the index of a recursive graph.
We can color this graph by mapping every node of $P$ to $0$, and mapping
all other nodes to their integer codes.  Since $0$ is used infinitely often
in this coloring, $f(e) \in \Cal G$.

Now suppose $e \notin \Cal T$.  If $e$ is not the index of a recursive tree,
then $f(e)$ is not the index of a recursive graph, so $f(e) \notin \Cal G$.
If we suppose that $e$ is the index of a recursive tree $T$, then $T$ is
well-founded.  Suppose, by way of contradiction, that there is a coloring
of the associated recursive graph $G_e$ that uses $0$ infinitely often.
All the nodes of $G_e$ that are colored $0$ correspond to comparable nodes
of $T$, contradicting the claim that $T$ is well-founded.  Again, we have
$f(e) \notin \Cal G$, completing the proof that $\Cal T \le_1 \Cal G$. \qed
\enddemo

The techniques of reverse mathematics can be used to draw a distinction
between the first two of our infinite analogs and the third.  The
following two results make use of the axiom systems
$\rca$ ({\sl R}ecursive {\sl C}omprehension {\sl A}xiom),
$\aca$ ({\sl A}rithmetical {\sl C}omprehension {\sl A}xiom),
and $\p11$ ($\Pi^1_1$ {\sl C}omprehension {\sl A}xiom).
For a brief overview of reverse mathematics, see {\smc Simpson} \cite{11}

%Insert Theorem D (now 5) and Theorem E (now 6) here.
\proclaim{Theorem 5 ($\rca$)}
The following are equivalent:
\roster
\item
$\aca$.
\item
For any sequence of graphs $\langle G_i : i \in \omega \rangle$,
there is a function $s:\omega\rightarrow 2$ such that $s(i) = 1$ if and only if
$G_i$ is 3-chromatic.
\item
For any sequence of graphs $\langle G_i : i \in \omega \rangle$,
there is a function $s:\omega\rightarrow 2$ such that $s(i) = 1$ if and only if
every connected component of $G_i$ is finitely colorable.
\endroster
\endproclaim

\demo{Proof}
To prove \therosteritem1$\rightarrow$\therosteritem2 and
\therosteritem1$\rightarrow$\therosteritem3, it suffices
to show that the function $s$ is arithmetically definable
in $\langle G_i : i \in \omega \rangle$.  For \therosteritem2,
imitating the first paragraph of the
proof of Theorem 1 yields a $\Pi^0_1$ defining
formula.  Similarly, for \therosteritem3, imitating the proof
of Theorem 2 yields a $\Pi^0_3$ defining formula.

By Lemma 2.7 of \cite{10},
to prove that \therosteritem2$\rightarrow$\therosteritem1
and \therosteritem3$\rightarrow$\therosteritem1, it suffices to show that
given any injection $g:\omega \rightarrow \omega$, $\rca$ can
prove the existence of a sequence of graphs
$\langle G_i : i \in \omega\rangle$ such that the range of $g$ is
$\Delta^0_1$ definable in the associated function $s$.  Fix $g$
and assume $\rca$.  We will define a sequence of graphs that works
for both \therosteritem2 and \therosteritem3.  Let $G_n$ have
$\omega$ as its vertex set.  For every $j \in \omega$,
include the edge $(j,j+1)$ in $G_n$.  For $j<k$, add the edge
$(j,k)$ to $G_n$ if and only if
$\exists t \le j \, ( g(t)=n)$.  The sequence
$\langle G_i : i \in \omega \rangle$ is
$\Delta^0_1$ definable in $g$, so $\rca$ proves it exists.
Let $s$ be as in \therosteritem2 or \therosteritem3.
Then $s(n)=1$ if and only if $n$ is not in the range of $g$.
Thus, the range of $g$ is $\Delta_1^0$ definable in $s$,
as desired.  \qed
\enddemo

\proclaim{Theorem 6 ($\rca$)}
The following are equivalent:
\roster
\item
$\p11$.
\item
For any sequence of graphs $\langle G_i : i \in \omega \rangle$,
there is a function $s:\omega\rightarrow 2$ such that $s(i) = 1$ if and only if
$G_i$ has a coloring in which one color is used infinitely often.
\endroster
\endproclaim

\demo{Proof}
To prove that \therosteritem1$\rightarrow$\therosteritem2, it suffices
to note that the function $s$ is $\Sigma_1^1$ definable in
$\langle G_i : i \in \omega \rangle$, and so exists by
$\p11$.  To prove the converse, we will use the fact
that $\p11$ is equivalent to the existence of a function that
decides which members of a sequence of trees are well founded.
(This is an easy consequence of Lemma 6.1 in \cite{3}.)
Assume
$\rca$, and suppose that $\langle T_i : i\in \omega \rangle$ is
a sequence of trees.  With each tree $T_n$, we associate a graph
$G_n$ as follows.  The vertices of $G_n$ are the nodes of $T_n$,
and two vertices of $G_n$ are connected if and only if the associated
nodes are incomparable in the tree ordering.  The sequence
$\langle G_i : i \in \omega \rangle$ is $\Delta^0_1$ definable in
$\langle T_i : i \in \omega \rangle$, and so exists by $\rca$.
Let $s$ be as in \therosteritem2.  Then $s(i)=1$ if and only if
$G_i$ contains an infinite collection of pairwise disconnected
vertices, which occurs if and only if $T_i$ is not well founded.
Thus \therosteritem2 implies $\p11$, completing the proof.  \qed
\enddemo

\subhead{Other Variability}\endsubhead

From the results in the preceding section, it is clear
that the recursion theoretic strength of infinite analogs
depends in part on their formulation.  As shown by
{\smc Harel} and {\smc Tirza Hirst} \cite{7}, adoption of a
standardized translation yields interesting parallels
between finite complexity and recursion theoretic complexity
for restricted classes of problems.  However, for broader
classes of problems, the parallels break down.  In this
section, we will consider three problems of diverse
finite complexity that all have $\Sigma^1_1$ complete
infinite analogs, thus illustrating our second theme.

Consider the following three variants of the subgraph
isomorphism problem:

\proclaim{P1}
Given a pair of finite graphs, $H$ and $G$, determine if $H$ is
isomorphic to a subgraph of $G$.
\endproclaim

\proclaim{P2}
For a fixed finite graph $H$, given a finite graph $G$, determine if $H$ is
isomorphic to a subgraph of $G$.
\endproclaim

\proclaim{P3}
For a fixed finite graph $G$, given a finite graph $H$, determine if $H$ is
isomorphic to a subgraph of $G$.
\endproclaim

{\bf P1} is the familiar form of the subgraph isomorphism problem,
and is known to be NP complete \cite{2}.
One algorithm for solving {\bf P2} and {\bf P3} consists of enumerating
all functions from $H$ into $G$, and checking each one to see if it is the
desired isomorphism.  The number of functions to check is bounded
by $|G|^{|H|}$, where $|G|$ denotes the number of vertices of $G$.  Since
$H$ is fixed in {\bf P2}, the number of functions to check is a constant
power of $|G|$.  Furthermore, the number of steps required to check each
function is bounded by a constant based on the fixed value $|H|$.
Thus, {\bf P2} can be solved in a number of steps which is
bounded by a polynomial
in $|G|$.  In {\bf P3}, $G$ is fixed, and we can discard any graphs
$H$ such that $|H|>|G|$, so the number of steps required to solve an
instance of {\bf P3} is bounded by a constant based on the fixed value $|G|$.
Summarizing, the complexity of three problems ranges from NP complete
to constant time computable.

Compared to the coloring problem in \S 1, these subgraph isomorphism
problems have very straightforward infinite analogs.  Despite the
variation in the computational complexity of the finite problems,
their infinite analogs are all $\Sigma^1_1$ complete, as is shown
in the following three theorems.

%???Insert theorems F (now 7) G (now 8) and H (now 9).

\proclaim{Theorem 7}
{\rm ({\smc Tirza Hirst} and {\smc Harel} \cite{7})}
The set of indices of ordered pairs of recursive graphs, $(H,G)$, such
that $H$ is isomorphic to a subgraph of $G$ is $\Sigma_1^1$ complete.
\endproclaim

\demo{Proof}
Let $\Cal G$ be the set of indices of ordered pairs of recursive graphs
such that the first graph is isomorphic to a subgraph of the second.
Since $x \in \Cal G$ if and only if an appropriate isomorphism exists,
it is easy to see that $\Cal G$ is $\Sigma_1^1$ definable.

To prove that $\Cal G$ is $\Sigma_1^1$ complete, we will show that
$\Cal T \le_1 \Cal G$, where $\Cal T$ denotes the set of indices
of recursive trees which are not well founded.  With each $e\in\omega$,
we associate a pair of partial recursive graphs, $H_e$ and $G_e$.
$H_e$ is a countably infinite linear graph with a triangle
attached at one end.  To be precise,
the vertex set of $H_e$ is $\{v_n : n \in \omega \}$ and the
edge set is
$\{(v_0 , v_2 )\} \cup \{ (v_n , v_{n+1} ) : n \in \omega \}$.
If $e$ is the index of a recursive tree $T$, then
$G_e$ consists of a copy of $T$ with a triangle attached to
the root, and a collection of disconnected vertices.  In general,
the vertex set for $G_e$ consists of $\{v_0 , v_1 , v_2 \}$ and
(codes for) the elements of $\omega^{<\omega}$.  Let
$\sigma_0$ denote the code for the empty sequence.  The edge $(v_0 , \sigma _0 )$
and the three edges of the form $( v_i , v_j )$ where $i \ne j$
are included in $G_e$.  For every $\sigma$ and $\tau$ in $\omega^{<\omega}$,
the edge $(\sigma , \tau )$ is included in $G_e$ if and only if
$$
\eta_e (\sigma ) = \eta_e (\tau ) = 1 \land
\sigma \subseteq \tau \land
\neg \exists \alpha ( \sigma \subsetneq \alpha \subsetneq \tau ),
$$
where $\eta_e$ is the function defined in \S 1.  By the
{\it s-m-n} theorem, there is a recursive 1-1 function $f$ such
that for every $e$, $f(e)$ is an index for the pair
$( H_e , G_e )$.

If $e \in \Cal T$, then $e$ is the index of a recursive tree containing
an infinite path $P$.  In this case, $H_e$ is isomorphic to the subgraph
of $G_e$ consisting of the base triangle and a copy of $P$.  Thus $f(e) \in \Cal G$.

Now suppose that $e \notin \Cal T$.  If $e$ is not the index of a recursive
tree, then $G_e$ is not a recursive graph, so $f_e \notin \Cal G$.
If $e$ is the index of a recursive tree $\Cal T$, then
$\Cal T$ is well founded.  The graph $G_e$ is a copy of $T$ with
a triangle attached to its base.  Any isomorphism mapping
$H_e$ into $G_e$ must map the triangle in $H_e$ into the
triangle in $G_e$, and the linear portion of $H_e$ to an infinite
path in the copy of $T$.  Since $T$ is well founded, no such
isomorphism exists.  Thus $f(e) \notin \Cal G$, completing the
proof that $\Cal T \le_1 \Cal G$.  \qed
\enddemo

\proclaim{Theorem 8}
There is a recursive graph $H$, such that the set of indices of recursive
graphs containing a subgraph isomorphic to $H$ is $\Sigma_1^1$ complete.
\endproclaim

\demo{Proof}
In the proof of Theorem 7, $H_e$ is a fixed recursive graph defined
without reference to $e$.  Any recursive 1-1 function mapping $e$ to
an index for the graph $G_e$ (defined as in the proof of Theorem 7)
witnesses the desired 1-reduction.  \qed
\enddemo

\proclaim{Theorem 9}
There is a recursive graph $G$, such that the set of indices of recursive
graphs that are isomorphic to a subgraph of $G$ is $\Sigma_1^1$ complete.
\endproclaim

\demo{Proof}
We begin the proof by constructing
the recursive graph $G$.
This graph will consist of a countable
collection of subgraphs $\langle G_e : e \in \omega \rangle$,
where each $G_e$ consists of a tree-like substructure together
with some spurious disconnected subgraphs.

For each $e \in \omega$, $G_e$ will be constructed from cycles
labeled $C ( e, \sigma , k )$ for each non-empty
$\sigma \in \omega^{<\omega}$ and each $k \in \omega$.
The cycle $C ( e, \sigma , k )$ consists of $2(e+1)+2$
vertices joined to make a circular graph.  We designate
two vertices of $C ( e, \sigma , k )$ as
$v^0_{e, \sigma , k}$ and $v^1_{e, \sigma , k}$, and require
that the paths joining them contain $e+2$ edges.  To give
a concrete example, $C ( 1, \sigma , k )$
looks like a hexagon, with the bottom vertex labeled
$v^0_{1, \sigma , k}$ and the top vertex labeled
$v^1_{1, \sigma , k}$.

The tree-like substructure of $G_e$ consists of a triangular
base with a vertex labeled $t_0$, and branches consisting
of linked cycles.  We say that a cycle $C ( e, \sigma , k )$
is {\sl exact} if $k$ is the least integer such that
1) $\eta_e (\tau ) \downarrow$ by stage $k$ for every $\tau$ which
is an initial subsequence of $\sigma$ or has a code less than $\sigma$,
and 2) $\eta_e (\sigma ) = 1$.  (Here $\eta_e$ is the function defined
in \S 1.)  Edges are added to $G_e$ by the following two
rules.  Connect $v^0_{e,\sigma , k}$ to $t_0$ if and only if
$C ( e, \sigma , k )$ is an exact cycle and $\sigma$ is a sequence
of length 1.  Connect $v^1_{e,\sigma , k}$ to $v^0_{e,\tau , j}$
if and only if $C ( e, \sigma , k )$ and $C ( e, \tau , j )$
are exact cycles and $\tau = \sigma * \langle m \rangle$ for some
$m \in \omega$.  Cycles which are not exact are {\sl spurious};
they are included in $G_e$, but are never connected to the
tree-like substructure.

Let $G$ be the union of all the $G_e$'s.  $G$ is recursive, since
the rules for adding edges involve only bounded computations.
Furthermore, if $e$ is the code of a recursive tree $T$,
then the tree-like substructure of $G_e$ can be mapped into
$T$ by identifying exact cycles with corresponding nodes.
Viewing the cycles as nodes, the substructure is well founded
if and only if $T$ is a well founded tree.  If $e$ is not the
code of a recursive tree, $\eta_e$ is not total, and the
tree-like substructure of $G_e$ is finite.

Let $\Cal G$ be the set of indices of recursive graphs that
are isomorphic to a subgraph of $G$.  Since $x \in \Cal G$
if and only if an isomorphism exists, it is easy to see that
$\Cal G$ is $\Sigma_1^1$ definable.  To prove that
$\Cal G$ is $\Sigma_1^1$ complete, we will show that
$\Cal T \le_1 \Cal G$, where $\Cal T$ denotes the
set of indices of recursive trees which are not well founded.
With each $e\in \omega$, we associate a recursive graph $H_e$
consisting of a countable linear graph with each node replaced by
a $2(e+1)+2$ cycle and with a triangle attached at one end.
More precisely, $H_e$ contains a triangle with one vertex labeled
$t_0$, and (copies of) the cycles $C(e,\langle 0 \rangle , k )$
for each $k \in \omega$.  To the edges already specified,
we add the edge $(t_0 , v^0_{e,\langle 0 \rangle , k})$ and
the edges $( v^1_{e,\langle 0 \rangle , k} , v^0_{e,\langle 0 \rangle , k+1} )$
for each $k \in \omega$.  By the {\it s-m-n} Theorem,
there is a recursive 1-1 function $f$ such that for
every $e$, $f(e)$ is an index for $H_e$.

If $e \in \Cal T$, then $e$ is the index of a recursive tree
containing an infinite path $P$.  In this case, $H_e$ is isomorphic
to the subgraph of $G_e$ consisting of the base triangle
and a copy of $P$ with nodes replaced by cycles.  Thus $f(e) \in \Cal G$.

Now suppose that $e \notin \Cal T$.  Note that because the size
of the cycles varies with $e$, if $H_e$ is isomorphic to a subgraph of
$G$, then $H_e$ is isomorphic to a subgraph of $G_e$.  Since
$e \notin \Cal T$, $G_e$ consists of disconnected cycles
and a well founded tree-like substructure.  If $H_e$ is isomorphic
to a subgraph of $G_e$, then the tree-like substructure of
$G_e$ contains an infinite path, yielding a contradiction.
Thus $f(e) \notin \Cal G$ completing the proof that
$\Cal T \le _1 \Cal G$.  \qed
\enddemo

Using the reverse mathematics framework, the preceding three theorems can
be lumped together into a single equivalence result.

%???Insert theorem J (now 10)

\proclaim{Theorem 10 ($\rca$)}
The following are equivalent:
\roster
\item
$\p11$.
\item
For any sequence of ordered pairs of graphs,
$\langle (H_i , G_i ) : i \in \omega \rangle$,
there is a function $s:\omega\rightarrow 2$ such that $s(i) =1$
if and only if $H_i$ is isomorphic to a subgraph of $G_i$.
\item
For any graph $H$, and any sequence of graphs
$\langle G_i : i \in \omega \rangle$,
there is a function $s:\omega\rightarrow 2$ such that $s(i) =1$
if and only if $H$ is isomorphic to a subgraph of $G_i$.
\item
For any graph $G$, and any sequence of graphs
$\langle H_i : i \in \omega \rangle$,
there is a function $s:\omega\rightarrow 2$ such that $s(i) =1$
if and only if $H_i$ is isomorphic to a subgraph of $G$.
\endroster
\endproclaim

\demo{Proof}
To prove that \therosteritem1 implies
\therosteritem2, \therosteritem3, or \therosteritem4,
it suffices to note that the function $s$ is
$\Sigma_1^1$ definable in the appropriate sequence of
graphs.  Since \therosteritem3 is a special case
of \therosteritem2, we need only show that
\therosteritem3$\rightarrow$\therosteritem1 and
\therosteritem4$\rightarrow$\therosteritem1 to
complete the proof.  As in the proof of Theorem 6,
we will determine which members of a sequence of trees
are well founded.  For the remainder of the proof,
assume $\rca$ and let $\langle T_i : i \in \omega \rangle$
be a sequence of trees.

To prove that \therosteritem3$\rightarrow$\therosteritem1,
we use a simplified version of the construction in the
proof of Theorem 8.  As in that proof, let $H$ be
a countable linear graph with a triangle attached to
one end.  For each $n \in \omega$, let
$G_n$ be a copy of $T_n$, with a triangle attached to
the root.  The graph $H$ and the sequence
$\langle G_i : i \in \omega \rangle$ are
$\Delta^0_1$ definable in
$\langle T_i : i \in \omega \rangle$, so
$\rca$ proves that they exist.  Let $s$ be as in
\therosteritem3.  Then $s(i)=1$ if and only if
$H$ is isomorphic to a subgraph of $G_i$, which
occurs if and only if $T_i$ has an infinite path.
Thus \therosteritem3 implies $\p11$.

To prove that \therosteritem4$\rightarrow$\therosteritem1,
we use a simplified version of the proof of Theorem 9.
As in that proof, let $H_n$ consist of a linear graph
with each node replaced by a $2(n+1)+2$ cycle, and with
a triangle attached to one end.  The graph $G$ consists
of subgraphs $G_n$ for each $n \in \omega$, where each $G_n$ is
a copy of $T_n$ with non-base nodes replaced by
$2(n+1)+2$ cycles, and a triangle attached to the base
node.  The graph $G$ and the sequence
$\langle H_i : i \in \omega \rangle$ are $\Delta_1^0$
definable in $\langle T_i : i \in \omega \rangle$, so
$\rca$ proves that they exist.  If $s$ is as in
\therosteritem4, then $s(i)=1$ if and only if $H_i$ is
isomorphic to a subgraph of $G$, which occurs if and
only if $H_i$ is isomorphic to a subgraph of $G_i$.
Finally, $H_i$ is isomorphic to a subgraph of $G_i$
if and only if $T_i$ is not well founded, so
\therosteritem4 implies $\p11$, completing the proof.  \qed
\enddemo

Although infinite analogs are useful for studying restricted classes
of problems, the preceding examples indicate that, in a general setting,
their behavior does not necessarily parallel that of the associated
finite problems.  However, examination of results in finite complexity
can provide motivation for appealing results in recursion theory and
reverse mathematics.

\enddocument
\end